\documentclass{amsart}
\usepackage{amsfonts}
\usepackage{amsmath}
\usepackage{amssymb}
\usepackage{amsthm}

\newtheorem{lemma}{Lemma}
\newtheorem{theorem}{Theorem}

\newtheorem{exmp}{Example}
\newtheorem{cor}{Corollary}

\newtheorem{fact}{Fact}

\begin{document}

\author{Mark  Pankov}
\title{On geometry of symplectic involutions}
\address{Institute of Mathematics NASU, Kiev}
\email{pankov@imath.kiev.ua}

\begin{abstract}
Let $V$ be a $2n$-dimensional vector space over a field $F$
and $\Omega$ be a non-degenerate symplectic form on $V$.
Denote by ${\mathfrak H}_{k}(\Omega)$
the set of all $2k$-dimensional subspaces $U\subset V$
such that the restriction $\Omega|_{U}$ is non-degenerate.
Our main result (Theorem 1) says that
if $n\ne 2k$ and $\max(k,n-k)\ge 5$ then any bijective transformation of
${\mathfrak H}_{k}(\Omega)$ preserving the class of base subsets
is induced by a semi-simplectic automorphism of $V$.
For the case when $n\ne 2k$ this fails,
but we have a weak version of this result (Theorem 2).
If the characteristic of $F$ is not equal to $2$ then
there is a one-to-one correspondence between elements of ${\mathfrak H}_{k}(\Omega)$
and symplectic $(2k,2n-2k)$-involutions
and Theorem 1 can be formulated as follows:
for the case when $n\ne 2k$ and $\max(k,n-k)\ge 5$
any commutativity preserving bijective transformation of
the set of symplectic $(2k,2n-2k)$-involutions can be extended
to an automorphism of the symplectic group.
\end{abstract}

\maketitle

\section{Introduction}

Let $W$ be an $n$-dimensional vector space over a division ring $R$ and $n\ge 3$.
We put ${\mathcal G}_{k}(W)$ for the Grassmannian of $k$-dimensional subspaces of $W$.
The projective space associated with $W$ will be denoted by ${\mathcal P}(W)$.

Let us consider the set ${\mathfrak G}_{k}(W)$ of all pairs
$$(S,U)\in {\mathcal G}_{k}(W)\times {\mathcal G}_{n-k}(W),$$
where $S+U=W$.
If $B$ is a base for ${\mathcal P}(W)$ then
the {\it base subset} of ${\mathfrak G}_{k}(W)$ associated with the base $B$
consists of all $(S,U)$ such that $S$ and $U$ are spanned by elements of $B$.
{\it If $n\ne 2k$ then
any bijective transformation of ${\mathfrak G}_{k}(W)$ preserving
the class of base subsets
is induced by a semi-linear isomorphism of $W$ to itself or to
the dual space $W^{*}$} (for $n=2k$ this fails, but some weak version
of this result holds true).
Using Mackey's ideas \cite{Mackey}
J. Dieudonn\'{e} \cite{D1} and C. E. Rickart  \cite{Rickart}
have proved this statement for $k=1,n-1$.
For the case when $1<k<n-1$ it was established by author \cite{Pankov}.
Note that adjacency preserving transformations of ${\mathfrak G}_{k}(W)$
were studied in \cite{HP}.

Now suppose that the characteristic of $R$ is not equal to $2$
and consider an involution $u\in{\rm GL}(W)$.
There exist two subspaces $S_{+}(u)$ and $S_{-}(u)$ such that
$$u(x)=x\;\mbox{ if }\;x\in S_{+}(u)\,,
\;\;\;u(x)=-x\;\mbox{ if }\;x\in S_{-}(u)$$
and
$$W=S_{+}(u)+S_{-}(u).$$
We say that $u$ is a $(k,n-k)$-{\it involution}
if the dimensions of $S_{+}(u)$ and $S_{-}(u)$ are equal to $k$ and $n-k$, respectively.
The set of $(k,n-k)$-involutions will be denoted by ${\mathfrak I}_{k}(W)$.
There is the natural one-to-one correspondence between
elements of ${\mathfrak I}_{k}(W)$ and ${\mathfrak G}_{k}(W)$
such that each base subsets of ${\mathfrak G}_{k}(W)$ corresponds to
a maximal set of mutually permutable $(k,n-k)$-involutions.
Thus any commutativity preserving transformation of ${\mathfrak I}_{k}(W)$
can be considered as a transformation of ${\mathfrak G}_{k}(W)$ preserving
the class of base subsets,
and our statement shows that {\it if $n\ne 2k$ then
any commutativity preserving bijective transformation of ${\mathfrak I}_{k}(W)$
can be extended to an automorphism of ${\rm GL}(W)$.}

In the present paper we give symplectic analogues of these results.

\section{Results}

\subsection{}
Let $V$ be a $2n$-dimensional vector space over a field $F$
and $\Omega:V\times V\to F$ be a non-degenerate symplectic form.
The form $\Omega$ defines on the set of subspaces of $V$ the orthogonal relation
which will be denoted by $\perp$.
For any subspace $S\subset V$ we put $S^{\perp}$ for the orthogonal complement to $S$.
A subspace $S\subset V$ is said to be {\it non-degenerate} if
the restriction $\Omega|_{S}$ is non-degenerate;
for this case $S$ is even-dimensional and $S+S^{\perp}=V$.
We put ${\mathfrak H}_{k}(\Omega)$ for the set of non-degenerate $2k$-dimensional subspaces.
Any element of ${\mathfrak H}_{k}(\Omega)$ can be presented as
the sum of $k$ mutually orthogonal elements of ${\mathfrak H}_{1}(\Omega)$.

Let us consider the projective space ${\mathcal P}(V)$ associated with $V$.
The points of this space are $1$-dimensional subspaces of $V$,
and each line consists of all
$1$-dimensional subspaces contained in a certain $2$-dimensional subspace.
A line of ${\mathcal P}(V)$ is called {\it hyperbolic}
if the corresponding $2$-dimensional subspace belongs to ${\mathfrak H}_{1}(\Omega)$;
otherwise, the line is said to be {\it isotropic}.
Points of ${\mathcal P}(V)$ together with the family of isotropic lines
form the well-known {\it polar space}.
Some results related with the {\it hyperbolic symplectic geometry}
(spanned by points of ${\mathcal P}(V)$ and hyperbolic lines)
can be found in \cite{Cuypers}, \cite{Gramlich}, \cite{Hall}.

A base $B=\{P_{1},\dots, P_{2n}\}$ of ${\mathcal P}(V)$ is called {\it symplectic} if
for any $i\in \{1,\dots,2n\}$ there is unique $\sigma(i)\in \{1,\dots,2n\}$
such that $P_{i}\not\perp P_{\sigma(i)}$.
Then the set ${\mathfrak S}_{1}$ consisting of all
$$S_{i}:=P_{i}+P_{\sigma(i)}$$
is said to be the {\it base subset} of ${\mathfrak H}_{1}(\Omega)$
associated with the base $B$.
For any $k\in \{2,\dots, n-1\}$ the set ${\mathfrak S}_{k}$ consisting of all
$S_{i_{1}}+\dots+S_{i_{k}}$
($i_{1},\dots, i_{k}$ are different)
will be called the {\it base subset} of ${\mathfrak H}_{k}(\Omega)$ associated with
${\mathfrak S}_{1}$ (or defined by ${\mathfrak S}_{1}$).

Now suppose that the characteristic of $F$ is not equal to $2$.
An involution $u\in {\rm GL}(V)$ is symplectic
(belongs to the group ${\rm Sp}(\Omega)$) if and only if
$S_{+}(u)$ and $S_{-}(u)$ are
non-degenerate and $S_{-}(u)=(S_{+}(u))^{\perp}$.
We denote by ${\mathfrak I}_{k}(\Omega)$
the set of symplectic $(2k, 2n-2k)$-involutions.
There is the natural bijection
$$i_{k}:{\mathfrak I}_{k}(\Omega)\to{\mathfrak H}_{k}(\Omega),\;\;\;\;\;
u\to S_{+}(u).$$
We say that ${\mathfrak X}\subset {\mathfrak I}_{k}(\Omega)$
is a $MC$-{\it subset}\/ if any two elements of ${\mathfrak X}$
are commutative and
for any $u\in {\mathfrak I}_{k}(\Omega)\setminus {\mathfrak X}$
there exists $s\in {\mathfrak X}$ such that $su\ne us$
(in other words, ${\mathfrak X}$ is a maximal set
of mutually permutable elements of ${\mathfrak I}_{k}(\Omega)$).

\begin{fact}\cite{D1}, \cite{D2}
${\mathfrak X}$ is a $MC$-subset of ${\mathfrak I}_{k}(\Omega)$
if and only if $i_{k}({\mathfrak X})$
is a base subset of ${\mathfrak H}_{k}(\Omega)$.
For any two commutative elements of ${\mathfrak I}_{k}(\Omega)$
there is a $MC$-subset containing them.
\end{fact}

Fact 1 shows that a bijective transformation $f$ of ${\mathfrak H}_{k}(\Omega)$
preserves the class of base subsets if and only if
$i^{-1}_{k}fi_{k}$ is commutativity preserving.

\subsection{}
If $l$ is an element of $\Gamma {\rm Sp}(\Omega)$
(the group of semi-linear  automorphisms preserving $\Omega$)
then for each number $k\in \{1,\dots,n-1\}$
we have the bijective transformation
$$(l)_{k}:{\mathfrak H}_{k}(\Omega)\to {\mathfrak H}_{k}(\Omega),\;\;\;\;\;
U\to l(U)$$
which preserves the class of base subsets.
The bijection
$$p_{k}:{\mathfrak H}_{k}(\Omega)\to {\mathfrak H}_{n-k}(\Omega),\;\;\;\;\;
U\to U^{\perp}$$
sends base subsets to base subsets.
We will need the following trivial fact.

\begin{fact}
Let $f$ be a bijective transformation of ${\mathfrak H}_{k}(\Omega)$
preserving the class of base subsets.
Then the same holds for the transformation $p_{k}fp_{n-k}$.
Moreover, if $f=(l)_{k}$ for certain $l\in \Gamma {\rm Sp}(\Omega)$
then $p_{k}fp_{n-k}=(l)_{n-k}$.
\end{fact}

Two distinct elements of  ${\mathfrak H}_{1}(\Omega)$ are orthogonal
if and only if there exists a base subset containing them,
thus for any bijective transformation $f$ of ${\mathfrak H}_{1}(\Omega)$
the following condition are equivalent:
\begin{enumerate}
\item[---] $f$ preserves the relation $\perp$,
\item[---] $f$ preserves the class of base subsets.
\end{enumerate}
It is not difficult to prove
(see \cite{D1}, p. 26-27 or \cite{Rickart}, p. 711-712)
that if one of these conditions  holds then
$f$ is induced by an element of $\Gamma {\rm Sp}(\Omega)$.
Fact 2 guarantees that the same is fulfilled for
bijective transformations of ${\mathfrak H}_{n-1}(\Omega)$ preserving the class of
base subsets.
This result was exploited by
J. Dieudonn\'{e} \cite{D1} and C. E. Rickart \cite{Rickart}
to determining automorphisms of the group ${\rm Sp}(\Omega)$.

\begin{theorem}
If $n\ne 2k$ and $\max(k,n-k)\ge 5$ then
any bijective transformation of ${\mathfrak H}_{k}(\Omega)$
preserving the class of base subsets is induced by
an element of $\Gamma {\rm Sp}(\Omega)$.
\end{theorem}

\begin{cor}
Suppose that the characteristic of $F$ is not equal to $2$.
If $n\ne 2k$ and $\max(k,n-k)\ge 5$ then any commutativity preserving
bijective transformation $f$ of ${\mathfrak I}_{k}(\Omega)$
can be extended to an automorphism of ${\rm Sp}(\Omega)$.
\end{cor}

\begin{proof}[Proof of Corollary]
By Fact 1, $i_{k}fi^{-1}_{k}$ preserves the class of base subsets.
Theorem 1 implies that $i_{k}fi^{-1}_{k}$ is induced by
$l\in \Gamma {\rm Sp}(\Omega)$.
The automorphism $u\to lul^{-1}$ is as required.
\end{proof}

\subsection{}
For the case when $n=2k$ Theorem 1 fails.

\begin{exmp}{\rm
Suppose that $n=2k$ and ${\mathfrak X}$ is a subset of
${\mathfrak H}_{k}(\Omega)$ such that for any $U\in {\mathfrak X}$
we have $U^{\perp}\in {\mathfrak X}$.
Consider the transformation of ${\mathfrak H}_{k}(\Omega)$
which sends each $U\in {\mathfrak X}$ to $U^{\perp}$
and leaves fixed all other elements.
This transformation preserves the class of base subsets
(any base subset of ${\mathfrak H}_{k}(\Omega)$ contains $U$ together with $U^{\perp}$),
but it is not induced by a semilinear automorphism if ${\mathfrak X}\ne\emptyset$.
}\end{exmp}

If $n=2k$ then we denote by $\overline{{\mathfrak H}}_{k}(\Omega)$
the set of all subsets
$\{U,U^{\perp}\}\subset {\mathfrak H}_{k}(\Omega)$.
Then every $l\in\Gamma {\rm Sp}(\Omega)$ induces the bijection
$$(l)'_{k}:\overline{{\mathfrak H}}_{k}(\Omega)\to \overline{{\mathfrak H}}_{k}(\Omega)
,\;\;\;\;\;
\{U,U^{\perp}\}\to \{l(U), l(U^{\perp})=l(U)^{\perp}\}.$$
The transformation from Example 1 gives the identical transformation
of $\overline{{\mathfrak H}}_{k}(\Omega)$.

\begin{theorem}
Let $n=2k\ge 14$ and $f$ be a bijective transformation of ${\mathfrak H}_{k}(\Omega)$
preserving the class of base subsets.
Then $f$ preserves the relation $\perp$
and induces a bijective transformation of
$\overline{{\mathfrak H}}_{k}(\Omega)$.
The latter mapping is induced by an element of $\Gamma {\rm Sp}(\Omega)$.
\end{theorem}

\begin{cor}
Let $n=2k\ge 14$ and
$f$ be a commutativity preserving bijective transformation of
${\mathfrak I}_{k}(\Omega)$.
Suppose also that the characteristic of $F$ is not equal to $2$.
Then there exists an automorphism $g$ of the group ${\rm Sp}(\Omega)$
such that $f(u)=\pm g(u)$ for any $u\in {\mathfrak I}_{k}(\Omega)$.
\end{cor}

\section{Inexact subsets}
In this section we suppose that $n\ge 4$ and $1<k<n-1$.

\subsection{Inexact subsets of ${\mathfrak G}_{k}(W)$}
Let $B=\{P_{1},\dots,P_{n}\}$ be a base of ${\mathcal P}(W)$.
For any $m\in \{1,\dots, n-1\}$ we denote by ${\mathfrak B}_{m}$
the base subset of ${\mathfrak G}_{m}(W)$ associated with $B$
(the definition was given in section 1).

If $\alpha=(M,N)\in {\mathfrak B}_{m}$ then we put
${\mathfrak B}_{k}(\alpha)$ for the set of all $(S,U)\in{\mathfrak B}_{k}$
where $S$ is incident to $M$ or $N$
(then $U$ is incident to $N$ or $M$, respectively),
the set of all $(S,U)\in{\mathfrak B}_{k}$
such that $S$ is incident to $M$ will be denoted by
${\mathfrak B}^{+}_{k}(\alpha)$.

A subset ${\mathfrak X}\subset {\mathfrak B}_{k}$
is called {\it exact} if there is only one base subset of ${\mathfrak G}_{k}(W)$
containing ${\mathfrak X}$; otherwise, ${\mathfrak X}$ is said to be
{\it inexact}.
If $\alpha \in {\mathfrak B}_{2}$
then ${\mathfrak B}_{k}(\alpha)$ is a maximal inexact subset of ${\mathfrak B}_{k}$
(Example 1 in \cite{Pankov}).
Conversely, we have the following.

\begin{lemma}[Lemma 2 of \cite{Pankov}]
If ${\mathfrak X}$ is a maximal inexact subset of ${\mathfrak B}_{k}$
then there exists $\alpha \in {\mathfrak B}_{2}$
such that ${\mathfrak X}={\mathfrak B}_{k}(\alpha)$.
\end{lemma}

\begin{lemma}[Lemmas 5 and 8 of \cite{Pankov}]
Let $g$ be a bijective transformation of ${\mathfrak B}_{k}$
preserving the class of maximal inexact subsets.
Then for any $\alpha \in {\mathfrak B}_{k-1}$ there exists
$\beta \in {\mathfrak B}_{k-1}$ such that
$$g({\mathfrak B}_{k}(\alpha))={\mathfrak B}_{k}(\beta);$$
moreover, we have
$$g({\mathfrak B}^{+}_{k}(\alpha))={\mathfrak B}^{+}_{k}(\beta)$$
if $n\ne 2k$.
\end{lemma}

\subsection{Inexact subsets of  ${\mathfrak H}_{k}(\Omega)$}
Let ${\mathfrak S}_1=\{S_{1},\dots,S_{n}\}$ be a base subsets of
${\mathfrak H}_{1}(\Omega)$.
For each number $m\in \{2,\dots, n-1\}$
we denote by ${\mathfrak S}_{m}$ the base subset of ${\mathfrak H}_{m}(\Omega)$
associated with ${\mathfrak S}_1$.

Let $M\in {\mathfrak S}_{m}$.
Then $M^{\perp}\in {\mathfrak S}_{n-m}$.
We put ${\mathfrak S}_{k}(M)$ for the set of all elements of
${\mathfrak S}_{k}$ incident to $M$ or $M^{\perp}$.
The set of all elements of ${\mathfrak S}_{k}$
incident to $M$ will be denoted by ${\mathfrak S}^{+}_{k}(M)$.

Let ${\mathfrak X}$ be a subset of ${\mathfrak S}_{k}$.
We say that ${\mathfrak X}$ is {\it exact} if
it is contained only in one base subset of ${\mathfrak H}_{k}(\Omega)$;
otherwise, ${\mathfrak X}$ will be called {\it inexact}.
For any $i\in \{1,\dots,n\}$ we denote by ${\mathfrak X}_{i}$ the set of all
elements of ${\mathfrak X}$ containing $S_i$.
If ${\mathfrak X}_{i}$ is not empty then we define
$$U_{i}({\mathfrak X}):=\bigcap_{U\in{\mathfrak X}_{i}}U,$$
and $U_{i}({\mathfrak X}):=\emptyset$ if ${\mathfrak X}_{i}$ is empty.
It is trivial that our subset is exact if
$U_{i}({\mathfrak X})=S_{i}$ for each $i$.

\begin{lemma}
${\mathfrak X}$ is exact if $U_{i}({\mathfrak X})\ne S_{i}$ only for one $i$.
\end{lemma}

\begin{proof}
Let ${\mathfrak S}'_1$ be a base subset of ${\mathfrak H}_{1}(\Omega)$
which defines a base subset of ${\mathfrak H}_{k}(\Omega)$ containing ${\mathfrak X}$.
If $j\ne i$ then $U_{j}({\mathfrak X})=S_{j}$
implies that $S_{j}$ belongs to ${\mathfrak S}'_1$.
Let us take $S'\in {\mathfrak S}'_1$ which does not coincide with any $S_{j}$, $j\ne i$.
Since $S'$ is orthogonal to all such $S_{j}$,
we have $S'=S_{i}$ and ${\mathfrak S}'_{1}={\mathfrak S}_{1}$.
\end{proof}

\begin{exmp}{\rm
Let $M\in {\mathfrak S}_{2}$. Then $M=S_{i}+S_{j}$ for some $i,j$.
We choose orthogonal $S'_{i},S'_{j}\in {\mathfrak H}_{1}(\Omega)$
such that $S'_{i}+S'_{j}=M$ and $\{S_{i},S_{j}\}\ne \{S'_{i},S'_{j}\}$.
Then
$$({\mathfrak S}_{1}\setminus \{S_{i},S_{j}\})\cup \{S'_{i},S'_{j}\}$$
is a base subset of ${\mathfrak H}_{1}(\Omega)$
which defines another base subset of ${\mathfrak H}_{k}(\Omega)$
containing ${\mathfrak S}_{k}(M)$.
Therefore, ${\mathfrak S}_{k}(M)$ is inexact.
Any $U\in {\mathfrak S}_{k}\setminus {\mathfrak S}_{k}(M)$
intersects $M$ by $S_{i}$ or $S_{j}$ and
$$U_{p}({\mathfrak S}_{k}(M)\cup \{U\})=S_{p}$$
if $p=i$ or $j$;
the same holds for all $p\ne i,j$.
By Lemma 3, ${\mathfrak S}_{k}(M)\cup \{U\}$ is exact
for any $U\in {\mathfrak S}_{k}\setminus {\mathfrak S}_{k}(M)$.
Thus the inexact subset ${\mathfrak S}_{k}(M)$ is maximal.
}\end{exmp}

\begin{lemma}
Let ${\mathfrak X}$ be a maximal inexact subset of ${\mathfrak S}_{k}$.
Then ${\mathfrak X}={\mathfrak S}_{k}(M)$
for certain $M \in {\mathfrak S}_{2}$.
\end{lemma}

\begin{proof}
By the definition,
there exists another base subset of ${\mathfrak H}_{k}(\Omega)$ containing ${\mathfrak X}$;
the associated base subset of ${\mathfrak H}_{1}(\Omega)$
will be denoted by ${\mathfrak S}'_{1}$.
Since our inexact subset is maximal,
we need to prove the existence of $M \in {\mathfrak S}_{2}$
such that ${\mathfrak X}\subset{\mathfrak S}_{k}(M)$.

Let us consider $i\in \{1,\dots,n\}$ such that $U_{i}$ is not empty
(from this moment we write $U_{i}$ in place of $U_{i}({\mathfrak X})$).
We say that the number $i$ is of {\it first type} if
the inclusion $U_{j}\subset U_{i}$, $j\ne i$ implies that
$U_{j}=\emptyset$ or $U_{j}=U_{i}$.
If $i$ is not of first type and the inclusion $U_{j}\subset U_{i}$, $j\ne i$
holds only for the case when $U_{j}=\emptyset$ or $j$ is of first type
then $i$ is said to be of {\it second type}.
Similarly, other types of numbers can be defined.

Suppose that there exists a number $j$ of first type such that $\dim U_{j}\ge 4$.
Then $U_{j}$ contains certain $M\in {\mathfrak S}_{2}$.
Since $j$ is of first type,
for any $U\in {\mathfrak X}$ one of the following possibilities is realized:
\begin{enumerate}
\item[---] $U\in {\mathfrak X}_{j}$ then $M\subset U_{j}\subset U$,
\item[---] $U\in {\mathfrak X}\setminus {\mathfrak X}_{j}$ then
$U\subset U^{\perp}_{j}\subset M^{\perp}$.
\end{enumerate}
This means that $M$ is as required.

Now suppose that $U_{j}=S_{j}$ for all $j$ of first type,
so $S_{j}\in {\mathfrak S}'_{1}$ if $j$ is of first type.
Consider any number $i$ of second type.
If $U_{i}\in {\mathfrak S}_m$ then $m\ge 2$ and there are exactly
$m-1$ distinct $j$ of first type such that $S_{j}=U_{j}$ is contained in $U_{i}$;
since all such $S_{j}$ belong to ${\mathfrak S}'_{1}$
and $U_{i}$ is spanned by elements of ${\mathfrak S}'_{1}$,
we have $S_{i}\in {\mathfrak S}'_{1}$.
Step by step we establish the same for other types.
Thus $S_{i}\in{\mathfrak S}'_{1}$ if $U_{i}$ is not empty.
Since ${\mathfrak X}$ is inexact,
Lemma 3 implies the the existence of two distinct numbers
$i$ and $j$ such that $U_{i}=U_{j}=\emptyset$.
We define $M:=S_{i}+S_{j}$.
Then any element of ${\mathfrak X}$ is contained in $M^{\perp}$
and we get the claim.
\end{proof}

Let ${\mathfrak S}'_{1}$ be another base subset of ${\mathfrak H}_{1}(\Omega)$
and ${\mathfrak S}'_{m}$, $m\in \{2,\dots,n-1\}$ be the base subset
of ${\mathfrak H}_{m}(\Omega)$ defined by ${\mathfrak S}'_{1}$.

\begin{lemma}
Let $h$ be a bijection of ${\mathfrak S}_{k}$
to ${\mathfrak S}'_{k}$ such that $h$ and $h^{-1}$
send maximal inexact subsets to maximal inexact subsets.
Then for any $M \in {\mathfrak S}_{k-1}$ there exists $M' \in {\mathfrak S}'_{k-1}$
such that
$$h({\mathfrak S}_{k}(M))={\mathfrak S}'_{k}(M');$$
moreover, we have
$$h({\mathfrak S}^{+}_{k}(M))={\mathfrak S}'^{+}_{k}(M')$$
if $n\ne 2k$.
\end{lemma}

\begin{proof}
Let ${\mathfrak B}_{m}$, $m\in \{1,\dots, n-1\}$ be as in subsection 3.1.
For each $m$ there is the natural bijection
$b_{m}:{\mathfrak B}_{m}\to {\mathfrak S}_{m}$
sending $(S,U)\in {\mathfrak B}_{m}$, $S=P_{i_{1}}+\dots +P_{i_{m}}$
to $S_{i_{1}}+\dots +S_{i_{m}}$.
For any $M\in {\mathfrak S}_{m}$ we have
$${\mathfrak S}_{k}(M)=b_{k}({\mathfrak B}_{k}(b^{-1}_{m}(M)))\;
\mbox{ and }\;
{\mathfrak S}^{+}_{k}(M)=b_{k}({\mathfrak B}^{+}_{k}(b^{-1}_{m}(M))).$$
Let $b'_{m}$ be the similar bijection of ${\mathfrak B}_{m}$ to ${\mathfrak S}'_{m}$.
Then $(b'_{k})^{-1}hb_{k}$ is a bijective transformation of ${\mathfrak B}_{k}$
preserving the class of base subsets
and our statement follows from Lemma 2.
\end{proof}

\section{Proof of Theorems 1 and 2}
By Fact 2, we need to prove Theorem 1 only for $k<n-k$.
Throughout the section we suppose that $1<k\le n-k$ and $n-k\ge 5$;
for the case when $n=2k$ we require that $n\ge 14$.

\subsection{}
Let $f$ be a bijective transformation of
${\mathfrak H}_{k}(\Omega)$ preserving the class of base subsets.
The restriction of $f$ to any base subset
satisfies the condition of Lemma 5.

For any subspace $T\subset V$ we denote by ${\mathfrak H}_{k}(T)$
the set of all elements of ${\mathfrak H}_{k}(\Omega)$ incident to $T$ or $T^{\perp}$,
the set of all elements of ${\mathfrak H}_{k}(\Omega)$ incident to $T$
will be denoted by ${\mathfrak H}^{+}_{k}(T)$.

In this subsection we show that Theorems 1 and 2 are
simple consequences of the following lemma.

\begin{lemma}
There exists a bijective transformation
$g$ of ${\mathfrak H}_{k-1}(\Omega)$ such that
$$g({\mathfrak H}^{+}_{k}(T))={\mathfrak H}^{+}_{k}(g(T))\;\;\;\;\;
\forall\; T\in {\mathfrak H}_{k-1}(\Omega)$$
if $n\ne 2k$, and
$$g({\mathfrak H}_{k}(T))=
{\mathfrak H}_{k}(g(T)))\;\;\;\;\;
\forall\; T\in {\mathfrak H}_{k-1}(\Omega)$$
for the case when $n=2k$.
\end{lemma}

Let ${\mathfrak S}_{k-1}$ be a base subset of ${\mathfrak H}_{k-1}(\Omega)$
and ${\mathfrak S}_{k}$ be the associated base subset of ${\mathfrak H}_{k}(\Omega)$
(these base subsets are defined by the same base subset of ${\mathfrak H}_{1}(\Omega)$).
By our hypothesis, $f({\mathfrak S}_{k})$ is a base subset of ${\mathfrak H}_{k}(\Omega)$;
we denote by ${\mathfrak S}'_{k-1}$ the associated base subset of
${\mathfrak H}_{k-1}(\Omega)$.
It is easy to see that $g({\mathfrak S}_{k-1})={\mathfrak S}'_{k-1}$,
so $g$  maps base subsets to base subsets.
Since $f^{-1}$ preserves the class of base subset, the same holds for $g^{-1}$.
Thus {\it $g$ preserves the class of base subsets}.

Now suppose that $g=(l)_{k-1}$ for certain $l\in \Gamma {\rm Sp}(\Omega)$.
Let $U$ be an element of ${\mathfrak H}_{k}(\Omega)$.
We take $M,N\in {\mathfrak H}_{k-1}(\Omega)$ such that
$U=M+N$.
If $n\ne 2k$ then
$$\{U\}={\mathfrak H}^{+}_{k}(M)\cap{\mathfrak H}^{+}_{k}(N)\;\mbox{ and }\;
\{f(U)\}={\mathfrak H}^{+}_{k}(l(M))\cap{\mathfrak H}^{+}_{k}(l(N)),$$
so $f(U)=l(M)+l(N)=l(U)$,
and we get $f=(l)_{k}$.
For the case when $n=2k$ we have
$$\{U, U^{\perp}\}=
{\mathfrak H}_{k}(M)\cap {\mathfrak H}_{k}(N)
\;\mbox{ and }\;
\{f(U), f(U)^{\perp}\}=
{\mathfrak H}_{k}(l(M))\cap {\mathfrak H}_{k}(l(N));$$
since $l(M)+l(N)=l(U)$ and
$l(M)^{\perp}\cap l(N)^{\perp}=(l(M)+l(N))^{\perp}=l(U)^{\perp},$
$$\{f(U), f(U)^{\perp}\}=\{l(U),l(U)^{\perp}\};$$
the latter means that $f=(l)'_{k}$.
Therefore, Theorem 1 can be proved by induction and
Theorem 2 follows from Theorem 1.

To prove Lemma 6 we use the following.

\begin{lemma}
Let $M\in {\mathfrak H}_{m}(\Omega)$ and $N$ be a subspace contained in $M$.
Then the following assertion are fulfilled:
\begin{enumerate}
\item[{\rm (1)}]
If $\dim N >m$ then $N$ contains an element of ${\mathfrak H}_{1}(\Omega)$.
\item[{\rm (2)}]
If $\dim N >m+2$ then $N$ contains two orthogonal elements of ${\mathfrak H}_{1}(\Omega)$.
\item[{\rm (3)}]
If $\dim N >m+4$ then $N$ contains three distinct mutually orthogonal elements of
${\mathfrak H}_{1}(\Omega)$.
\end{enumerate}
\end{lemma}

\begin{proof}
The form $\Omega|_{M}$ is non-degenerate.
If $\dim N >m$ then the restriction of $\Omega|_{M}$ to $N$ is non-zero.
This implies the existence of $S\in {\mathfrak H}_{1}(\Omega)$
contained in $N$.
We have
$$\dim N\cap S^{\perp}\ge \dim N-2,$$
and for the case when $\dim N >m+2$
there is an element of ${\mathfrak H}_{1}(\Omega)$
contained in $N\cap S^{\perp}$.
Similarly, (3) follows from (2).
\end{proof}

\subsection{Proof of Lemma 6 for $k<n-k$}
Let $T\in {\mathfrak H}_{k-1}(\Omega)$ and
${\mathfrak S}_{1}=\{S_{1},\dots,S_{n}\}$ be a base subset of ${\mathfrak H}_{1}(\Omega)$
such that
$$T^{\perp}=S_{1}+\dots+S_{n-k+1}\;\mbox{ and }\;T=S_{n-k+2}+\dots+S_{n}.$$
We put ${\mathfrak S}_{k}$ for the base subset of ${\mathfrak H}_{k}(\Omega)$
associated with ${\mathfrak S}_{1}$.
Then ${\mathfrak S}^{+}_{k}(T)$ consists of all
$$U_{i}:=T+S_{i},$$
where $i\in \{1,\dots, n-k+1\}$.
By Lemma 5, there exists $T'\in {\mathfrak H}_{k-1}(\Omega)$ such that
$$f({\mathfrak S}^{+}_{k}(T))\subset {\mathfrak H}^{+}_{k}(T').$$
We need to show that $f({\mathfrak H}^{+}_{k}(T))$ coincides with ${\mathfrak H}^{+}_{k}(T')$.

\begin{lemma}
Let $U\in {\mathfrak H}^{+}_{k}(T)$.
Suppose that there exist two distinct
$M,N\in {\mathfrak H}^{+}_{k}(T)$ such that
$f(M),f(N)$ belong to ${\mathfrak H}^{+}_{k}(T')$ and
there is a base subset of ${\mathfrak H}_{k}(\Omega)$
containing $M,N$ and $U$.
Then $f(U)$ is an element of ${\mathfrak H}^{+}_{k}(T')$.
\end{lemma}

\begin{proof}
If there exists a base subset of ${\mathfrak H}_{k}(\Omega)$ containing $M,N$ and $U$
then $T$ belongs to the associated base subset of ${\mathfrak H}_{k-1}(\Omega)$
and Lemma 5 implies the existence of $T''\in {\mathfrak H}_{k-1}(\Omega)$
such that $f(M),f(N)$ and $f(U)$ belong to ${\mathfrak H}^{+}_{k}(T'')$.
On the other hand, $f(M)$ and $f(N)$ are different elements of ${\mathfrak H}^{+}_{k}(T')$
and $f(M)\cap f(N)$ coincides with $T'$.
Hence $T'=T''$.
\end{proof}

For any $U\in {\mathfrak H}^{+}_{k}(T)$
we denote by $S(U)$ the intersection of $U$ and $T^{\perp}$,
it is clear that $S(U)$ is an element of ${\mathfrak H}_{1}(\Omega)$.

If $S(U)$ is contained in $S_{1}+\dots+S_{n-k-1}$
then $S(U),S_{n-k},S_{n-k+1}$ are mutually orthogonal
and there exists a base subset of ${\mathfrak H}_{k}(\Omega)$ containing
$U,U_{n-k},U_{n-k+1}$.
All $f(U_{i})$ belong to ${\mathfrak H}^{+}_{k}(T')$
and Lemma 8 shows that $f(U)\in {\mathfrak H}^{+}_{k}(T')$.

Let $U$ be an element of ${\mathfrak H}^{+}_{k}(T)$ such that
$S(U)$ is contained in $S_{1}+\dots+S_{n-k}$.
We have
$$\dim(S_{1}+\dots+S_{n-k-1})\cap S(U)^{\perp}\ge 2(n-k-2)>n-k-1$$
(the latter inequality follows from the condition $n-k\ge 5$)
and Lemma 7 implies the existence of $S'\in {\mathfrak H}_{1}(\Omega)$ contained in
$$(S_{1}+\dots+S_{n-k-1})\cap S(U)^{\perp}.$$
Then $S(U),S',S_{n-k+1}$ are mutually orthogonal and
there exists a base subset of ${\mathfrak H}_{k}(\Omega)$ containing
$U,T+S',U_{n-k+1}$.
It was proved above that $f(T+S')$ belongs to ${\mathfrak H}^{+}_{k}(T')$.
Since $f(U_{i})\in {\mathfrak H}^{+}_{k}(T')$ for each $i$,
Lemma 8 guarantees that $f(U)$ is an element of ${\mathfrak H}^{+}_{k}(T')$.

Now suppose that $S(U)$ is not contained in $S_{1}+\dots+S_{n-k}$.
Since $n-k\ge 5$,
$$\dim(S_{1}+\dots+S_{n-k})\cap S(U)^{\perp} \ge 2(n-k-1)>n-k+2.$$
By Lemma 7,
there exist two orthogonal $S',S''\in {\mathfrak H}_{1}(\Omega)$ contained in
$$(S_{1}+\dots+S_{n-k})\cap S(U)^{\perp}.$$
Then $S',S'',S(U)$ are mutually orthogonal and
there exists a base subset of ${\mathfrak H}_{k}(\Omega)$ containing
$S'+T, S''+T$ and $U$.
We have shown above that
$f(S'+T)$, $f(S''+T)$ belong to ${\mathfrak H}^{+}_{k}(T')$
and Lemma 8 shows that the same holds for $f(U)$.

So $f({\mathfrak H}^{+}_{k}(T))\subset {\mathfrak H}^{+}_{k}(T')$.
Since $f^{-1}$ preserves the class of base subsets,
the inverse inclusion holds true.
We define
$g:{\mathfrak H}_{k-1}(\Omega)\to {\mathfrak H}_{k-1}(\Omega)$
by $g(T):=T'$.
This transformation is bijective (otherwise, $f$ is not bijective).

\subsection{Proof of Lemma 6 for $n=2k$}
We start with the following.

\begin{lemma}
If $n=2k$ then $f(U^{\perp})=f(U)^{\perp}$ for any $U\in {\mathfrak H}_{k}(\Omega)$.
\end{lemma}

\begin{proof}
We take a base subset ${\mathfrak S}_{k}$ containing $U$.
Then $U^{\perp}\in {\mathfrak S}_{k}$.
Denote by ${\mathfrak S}_{k-1}$ the base subset of ${\mathfrak H}_{k-1}(\Omega)$
associated with ${\mathfrak S}_{k}$.
Let ${\mathfrak S}'_{k-1}$ be the base subset of ${\mathfrak H}_{k-1}(\Omega)$
associated with ${\mathfrak S}'_{k}:=f({\mathfrak S}_{k})$.
We choose $M,N\in {\mathfrak S}_{k-1}$ such that $U=M+N$.
Then
$$\{U,U^{\perp}\}={\mathfrak S}_{k}(M)\cap{\mathfrak S}_{k}(N)$$
and Lemma 5 guarantees that
$$\{f(U),f(U^{\perp})\}={\mathfrak S}'_{k}(M')\cap{\mathfrak S}'_{k}(N')$$
for some $M',N'\in {\mathfrak S}'_{k-1}$.
The set ${\mathfrak S}'_{k}(M')\cap{\mathfrak S}'_{k}(N')$
is not empty if one of the following possibilities is realized:
\begin{enumerate}
\item[---] $M'+N'$ and $M'^{\perp}\cap N'^{\perp}$ are elements of
${\mathfrak H}_{k-1}(\Omega)$
and ${\mathfrak S}'_{k}(M')\cap{\mathfrak S}'_{k}(N')$ consists
of these two elements.
\item[---]
$M'\subset N'^{\perp}$ and $N'\subset M'^{\perp}$,
then ${\mathfrak S}'_{k}(M')\cap{\mathfrak S}'_{k}(N')$ consists of
$4$ elements.
\end{enumerate}
Thus
$$\{f(U),f(U^{\perp})\}=\{M'+N',M'^{\perp}\cap N'^{\perp}\}.$$
Since $M'+N'$ and $M'^{\perp}\cap N'^{\perp}$ are orthogonal, we get the claim.
\end{proof}

Let $T\in {\mathfrak H}_{k-1}(\Omega)$.
As in the previous subsection we consider a base subset
${\mathfrak S}_{1}=\{S_{1},\dots,S_{n}\}$ of ${\mathfrak H}_{1}(\Omega)$
such that
$$T^{\perp}=S_{1}+\dots+S_{n-k+1}\;\mbox{ and }\;T=S_{n-k+2}+\dots+S_{n}.$$
We denote by ${\mathfrak S}_{k}$ the base subset of ${\mathfrak H}_{k}(\Omega)$
associated with ${\mathfrak S}_{1}$.
Then ${\mathfrak S}_{k}(T)$ consists of
$$U_{i}:=T+S_{i},\;\;\;i\in \{1,\dots, n-k+1\}$$
and their orthogonal complements.
Lemma 5 implies the existence of $T'\in {\mathfrak H}_{k-1}(\Omega)$
such that
$$f({\mathfrak S}_{k}(T))\subset {\mathfrak H}_{k}(T').$$
We show that $f(U)$ belongs to ${\mathfrak H}_{k}(T')$
for any $U\in {\mathfrak H}_{k}(T)$.

We need to establish this fact only for the case when
$U$ is an element of ${\mathfrak H}^{+}_{k}(T)$.
Indeed, if $U\in {\mathfrak H}^{+}_{k}(T^{\perp})$
then $U^{\perp}$ is an element of ${\mathfrak H}^{+}_{k}(T)$
and $f(U^{\perp})\in {\mathfrak H}_{k}(T')$ implies that
$f(U)=f(U^{\perp})^{\perp}$ belongs to ${\mathfrak H}_{k}(T')$.

\begin{lemma}
Let $U\in {\mathfrak H}^{+}_{k}(T)$.
Suppose that there exist distinct $M_{i}\in {\mathfrak H}^{+}_{k}(T)$, $i=1,2,3$
such that each $f(M_{i})$ belongs to ${\mathfrak H}_{k}(T')$ and
there is a base subset of ${\mathfrak H}_{k}(\Omega)$
containing $M_{1},M_{2},M_{3}$ and $U$.
Then $f(U)\in {\mathfrak H}_{k}(T')$.
\end{lemma}

\begin{proof}
By Lemma 5, there exists $T''\in {\mathfrak H}_{k-1}(\Omega)$ such that
$f(U)$, all $f(M_{i})$, and their orthogonal complements belong to ${\mathfrak H}_{k}(T'')$.
For any $i=1,2,3$ one of the subspaces $f(M_{i})$ or $f(M_{i})^{\perp}$
is an element of ${\mathfrak H}^{+}_{k}(T'')$;
we denote this subspace by $M'_{i}$.
Then
$$T''=\bigcap^{3}_{i=1}M'_{i}\;\;\mbox{ and }\;\;
T''^{\perp}=M'^{\perp}_{i}+M'^{\perp}_{j},\;i\ne j;$$
note also that the intersection of any $M'_{i}$ and $M'^{\perp}_{j}$
does not belong to ${\mathfrak H}_{k-1}(\Omega)$.
Since all $M'_{i}$ and $M'^{\perp}_{i}$ belong to
${\mathfrak H}_{k}(T')$, we have $T'=T''$.
\end{proof}

As in the previous subsection
for any $U\in {\mathfrak H}^{+}_{k}(T)$
we denote by $S(U)$ the intersection of $U$ and $T^{\perp}$,
it is an element of ${\mathfrak H}_{1}(\Omega)$.

If $S(U)$ is contained in $S_{1}+\dots+S_{n-k-2}$
then
$S(U),S_{n-k-1},S_{n-k},S_{n-k+1}$
are mutually orthogonal and
there exists a base subset of ${\mathfrak H}_{k}(\Omega)$ containing
$U,U_{n-k-1}$, $U_{n-k},U_{n-k+1}$.
Since $f(U_{i})\in {\mathfrak H}_{k}(T')$ for each $i$,
Lemma 10 shows that $f(U)$ belongs to ${\mathfrak H}_{k}(T')$.

Suppose that $S(U)$ is contained in $S_{1}+\dots+S_{n-k-1}$.
We have
$$\dim(S_{1}+\dots+S_{n-k-2})\cap S(U)^{\perp}\ge 2(n-k-3)>n-k-2$$
(since $k=n-k\ge 7$)
and Lemma 7 implies the existence of $S'\in {\mathfrak H}_{1}(\Omega)$ contained in
$$(S_{1}+\dots+S_{n-k-2})\cap S(U)^{\perp}.$$
Then $S(U),S',S_{n-k},S_{n-k+1}$ are mutually orthogonal,
so $U,T+S',U_{n-k},U_{n-k+1}$
are contained in a certain base subsets of ${\mathfrak H}_{k}(\Omega)$.
It was shown above that $f(T+S')$ is an element of ${\mathfrak H}_{k}(T')$
and Lemma 10 guarantees that $f(U)\in{\mathfrak H}_{k}(T')$
(recall that all $f(U_{i})$ belong to ${\mathfrak H}_{k}(T')$).

Consider the case when $S(U)$ is contained in $S_{1}+\dots+S_{n-k}$.
We have
$$\dim(S_{1}+\dots+S_{n-k-1})\cap S(U)^{\perp}\ge 2(n-k-2)>(n-k-1)+2$$
(recall that $k=n-k\ge 7$) and
there exist two orthogonal $S',S''\in {\mathfrak H}_{1}(\Omega)$ contained in
$$S_{1}+\dots+S_{n-k-1})\cap S(U)^{\perp}$$
(Lemma 7).
Then $S(U),S',S'',S_{n-k+1}$
are mutually orthogonal and there exists a base subsets of ${\mathfrak H}_{k}(\Omega)$
containing $U,T+S',T+S'',U_{n-k+1}$.
It follows from Lemma 10 that $f(U)\in {\mathfrak H}_{k}(T')$
(since $f(T+S')$, $f(T+S'')$ and any $f(U_{i})$ belong to ${\mathfrak H}_{k}(T')$).

Let $U$ be an element of ${\mathfrak H}_{k}(T')$
such that $S(U)$ is not contained in $S_{1}+\dots+S_{n-k}$.
Since $n=2k\ge 14$,
$$\dim(S_{1}+\dots+S_{n-k})\cap S(U)^{\perp}\ge 2(n-k-1)>n-k+4.$$
By Lemma 7, there exist mutually orthogonal
$S',S'',S'''\in {\mathfrak H}_{1}(\Omega)$ contained in
$$(S_{1}+\dots+S_{n-k})\cap S(U)^{\perp}.$$
A base subsets of ${\mathfrak H}_{k}(\Omega)$ containing
$U,T+S',T+S'',T+S'''$ exists.
It was shown above that
$f(T+S')$, $f(T+S'')$ and $f(T+S''')$ belong to ${\mathfrak H}_{k}(T')$
and Lemma 10 implies that the same holds for $f(U)$.

Thus $f({\mathfrak H}_{k}(T))\subset {\mathfrak H}_{k}(T')$.
As in the previous subsection we have the inverse inclusion
and define
$g:{\mathfrak H}_{k-1}(\Omega)\to {\mathfrak H}_{k-1}(\Omega)$
by $g(T):=T'$.

\end{document}